\documentclass[12pt]{elsarticle}

\usepackage{latexsym,amsfonts,amsmath,theorem,amssymb}
\usepackage{graphicx}

\def\be{\begin{equation}}
\def\bea{\begin{eqnarray*}}
\def\ee{\end{equation}}
\def\eea{\end{eqnarray*}}
\def\ba{\begin{array}}
\def\ea{\end{array}}
\def\bi{\begin{itemize}}
\def\ei{\end{itemize}}
\newtheorem{theo}{Theorem}
\newtheorem{lem}{Lemma}

\def\ZZ{{{\rm Z}\kern-.4em{\rm Z}}}
\def\RR{{{\rm I}\kern-.2em{\rm R}}}
\def\NN{{{\rm I}\kern-.2em{\rm N}}}
\def\TT{{{\rm T}\kern-.5em{\rm T}}}
\def\CC{{{\rm I}\kern-.5em{\rm C}}}

\journal{XXX}

\begin{document}
\begin{frontmatter}

\title{Nonconforming P1 elements on distorted triangulations: Lower bounds for the discrete energy norm error}
\author[bonn]{Peter Oswald}
\ead{agp.oswald@gmail.com}

\address[bonn]{Institute for Numerical Simulation (INS),
University of Bonn,
Wegelerstr. 6-8,
D-53115 Bonn
}
\date{}
\begin{abstract}
Compared to conforming P1 finite elements, nonconforming P1 finite element discretizations are thought to be less sensitive to the appearance of
distorted triangulations. E.g., optimal-order discrete $H^1$ norm best approximation error estimates for $H^2$ functions hold for arbitrary triangulations. However, the constants in similar estimates for the error of the Galerkin projection for second-order elliptic problems show a dependence on the maximum angle of all triangles in the triangulation. We demonstrate on the example of a special family of distorted triangulations that this dependence is essential, and due to the deterioration of the consistency error. We also provide examples of sequences of triangulations such that the nonconforming P1 Galerkin projections for a Poisson problem with polynomial solution do not converge or converge at arbitrarily low speed. The results complement analogous findings for conforming P1 finite elements.
\end{abstract}

\end{frontmatter}

\section{Introduction}\label{sec1}
Convergence estimates for the finite element method (FEM) in two and higher dimensions involve some shape regularity assumptions
for the underlying partitions. In two dimensions, to obtain optimal-order convergence estimates in the energy norm for triangular elements
when the maximal element diameter $h$ tends to zero,  the maximum angle condition introduced in \cite{BA,Ja} is sufficient. The natural question if this condition is also necessary has attracted less attention, even though
mesh generation strategies for the resolution of boundary and interior layers or discretizations involving moving meshes may lead to severely distorted triangle shapes. For conforming triangular P1 finite elements and the Poisson equation
\begin{equation}\label{Poisson}
-\Delta u =f,\qquad u\in H^1_0(\Omega),
\end{equation}
in \cite[Section 3]{BA} it was already shown on a particular example that the optimal-order O$(h)$ energy norm error estimate for smooth solutions $u\in H^2(\Omega)$ may not hold if the underlying sequence of triangulations severely violates the maximum angle condition. However, as was demonstrated in \cite{HKK}, there are many types of distorted triangulations violating the maximum angle condition but still admitting optimal-order error bounds for the Galerkin finite element method. In recent work \cite{Ku,Os}, some more precise statements about the necessity of the maximum angle condition for conforming  triangular P1 finite element  
discretizations have been made. E.g., in \cite{Os} for a particular Poisson problem on a square with polynomial solution, and a family of uniformly distorted triangulations already used in \cite{BA} and originating from \cite{Sc}, matching lower and upper bounds for the Galerkin energy norm error (or, equivalently, the error of best approximation by conforming P1 elements in the $H^1$ norm) have been obtained. These bounds precisely quantify the effect of the violation of the 
maximum angle condition on the convergence speed, and provide examples of sequences of triangulations where the Galerkin method does not converge to the
solution at all as $h\to 0$. In \cite{Ku}, a larger class of triangulations violating the maximum angle condition was investigated.

\begin{figure}[htb] 
\includegraphics[width=1.1\textwidth]{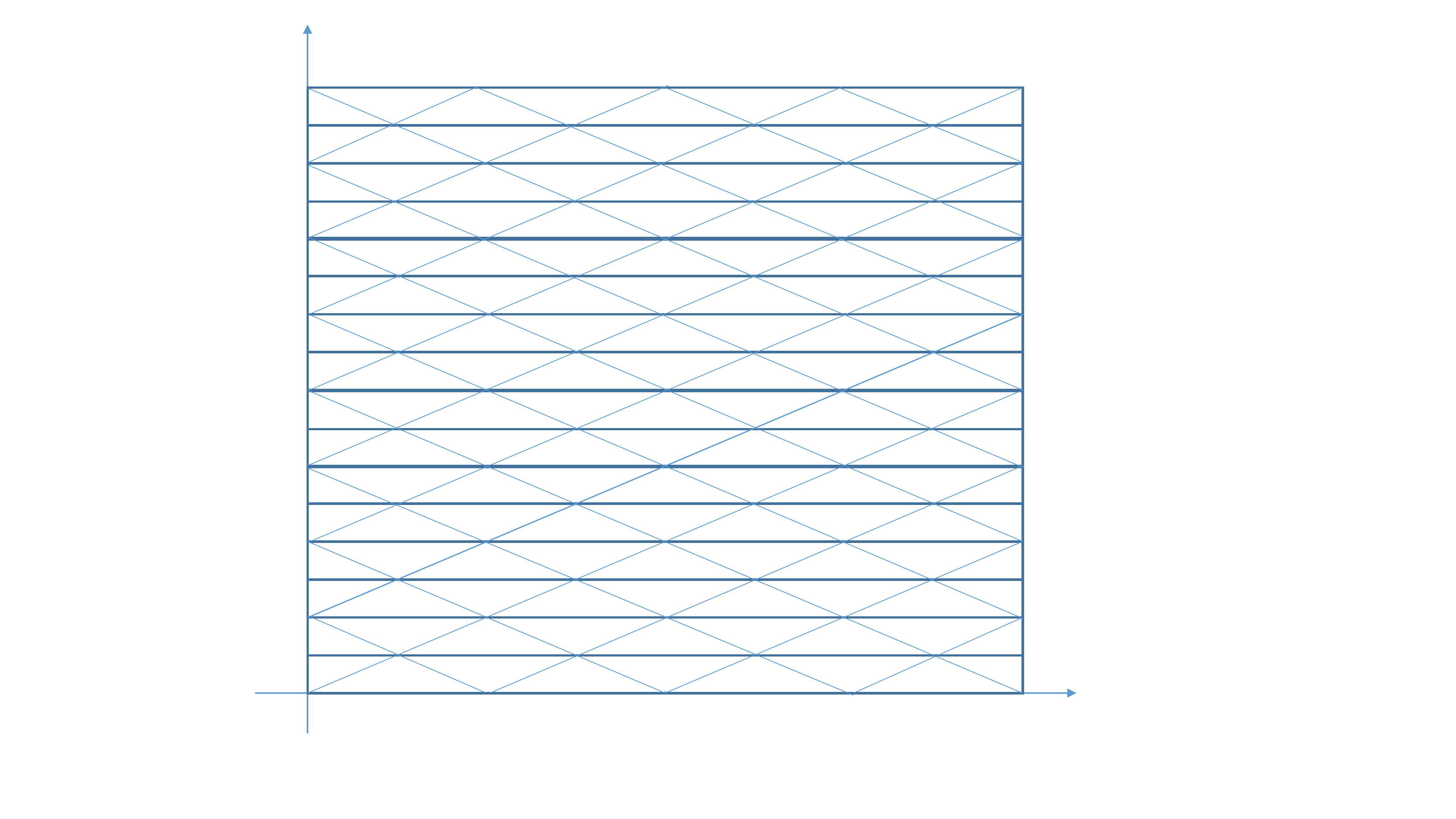}
\caption{Babu\v{s}ka-Aziz triangulation $\mathcal{T}_{4,8}$}
\label{fig1}
\end{figure}
One may wonder if the effects observed for conforming P1 elements in \cite{BA,HKK,Ku,Os} also hold for nonconforming and mixed finite element discretizations, where the maximum angle condition also figures as a sufficient condition, see \cite{AD,Br1,CGR} for a discussion of the lowest order Crouzeix-Raviart element \cite{CR} (commonly called nonconforming P1 element) and the closely related lowest order Raviart-Thomas element \cite{RT}. 
The examples from \cite{Os} show that the conforming P1 method may converge very slowly or even not converge at all while the nonconforming P1 method converges at the optimal O$(h)$ rate for $H^2$ solutions. From an approximation-theoretic point of view, the triangular nonconforming P1 element spaces generally provide better approximation properties in the discrete $H^1$ norm, independently of the shapes of triangles \cite{AD}, and the hope is that this may also extend to the error analysis of the Galerkin projection. However, since the Galerkin error of a nonconforming method also involves a 
consistency error, such an extension is not automatic, and, as it turns out, not possible. In the present paper, we modify the approach taken in \cite{Os},
and show for the same family $\{\mathcal{T}_{n,m}\}$ of triangulations as in \cite{BA,Os} that for the problem (\ref{Poisson}) with polynomial solution
$u(x,y)=x(1-x)y(1-y)$ and right-hand side $f(x,y)=2(x(1-x)+y(1-y))$ the nonconforming P1 Galerkin projections $u_{n,m}$ w.r.t. the triangulations $\mathcal{T}_{n,m}$
satisfy
\begin{equation}\label{TwoSided}
|u-u_{n,m}|_{H^1,\mathcal{T}_{n,m}} \approx \min(1,m/n^2), \qquad m\ge n>1.
\end{equation}
Here, $|\cdot|_{H^1,\mathcal{T}}$ stands for the discrete (sometimes called broken) $H^1$ norm associated with the triangulation $\mathcal{T}$, see Section \ref{sec2} for the definition. 
For $n=4$, $m=8$, the triangulation $\mathcal{T}_{n,m}$ is depicted in Figure \ref{fig1}. 
Since for $\mathcal{T}_{n,m}$ the mesh-size parameter $h$ equals $1/n$, and
the growth of $m/n$ measures the amount of deterioration of the maximum angle condition, we see that in general a violation of the maximum angle condition
immediately leads to a loss of convergence speed, and eventually to the loss of convergence, unless $m/n^2\to 0$ as $n\to \infty$.
However, examples in the spirit of \cite{HKK} show that not every sequence of triangulations containing irregularly shaped triangles share this behavior, and that the family $\{\mathcal{T}_{n,m}\}$ provides an extreme test case for the investigation of convergence problems with respect to distorted triangulations also
in the nonconforming P1 element case.

The two-sided estimate (\ref{TwoSided}) formally looks the same as the corresponding result from \cite{Os} for the conforming P1 element case but is different
in several aspects. First of all, the result from \cite{Os} is about the deterioration of the error of best approximation w.r.t. the conforming P1 element space on $\mathcal{T}_{n,m}$ for a Poisson problem with slightly different boundary conditions and with the polynomial solution $u(x,y)=x(1-x)/2$ depending only on the variable $x$. It can be checked that for problems with smooth solutions depending only on the variable $x$ the nonconforming P1 Galerkin projections for the triangulations on $\mathcal{T}_{n,m}$ converge at optimal speed O$(n^{-1})$, independently of the mesh distortion given by $m/n$ ($m\ge n$). We sketch the argument in Section \ref{sec4}. Thus, we need a truly two-dimensional approach. Secondly, the statement of (\ref{TwoSided}) is essentially about the consistency error induced by the nonconforming P1 element space on $\mathcal{T}_{n,m}$, and not about the best approximation error in the discrete $H^1$ norm.

The remainder of the paper is organized as follows. Section \ref{sec2} introduces notation and reviews the known upper estimates. In Section \ref{sec3}
the main result, the lower bound in (\ref{TwoSided}), is proved, some technical parts of this proof are delayed into appendices. The final Section \ref{sec4} offers complementary numerical evidence and contains some further remarks.

\section{Notation and Known Facts}\label{sec2}
Throughout the paper, we consider smooth solutions $u\in H^2(\Omega)\cap H_0^1(\Omega)$ of the Poisson problem (\ref{Poisson}) for a bounded polygonal domain $\Omega\subset \mathbb{R}^2$. Consequently, $f\in L_2(\Omega)$. Let $\mathcal{T}$ denote an arbitrary finite triangulation of $\Omega$ identified with a collection of closed triangles partitioning $\Omega$ with no hanging nodes. I.e., the intersection of any two  triangles in $\mathcal{T}$ is either empty or belongs to the vertex set $\mathcal{V}$ or to the edge set $\mathcal{E}$ of the triangulation.  Two characteristics of $\mathcal{T}$ are of interest to us: The mesh-width 
$$
h_{\mathcal{T}} := \max_{\Delta\in \mathcal{T}} h_\Delta,
$$
and the maximum angle
$$
\alpha_{\mathcal{T}}:= \max_{\Delta\in \mathcal{T}} \alpha_\Delta,
$$
where $h_\Delta$ denotes the length of the longest edge and $\alpha_\Delta$ the largest interior angle in a triangle $\Delta\in\mathcal{T}$, respectively.

The space of nonconforming P1 elements on $\mathcal{T}$ associated with homogeneous Dirichlet boundary conditions is denoted by $V_\mathcal{T}$,
and consists of all piecewise linear functions that are continuous across the midpoints of interior edges, and are zero at the midpoints of boundary edges. 
I.e., if $e\in \mathcal{E}$ is an interior edge shared by the  triangles  $\Delta^+$ and $\Delta^-$, then the two functions $v^\pm=v|
_{\Delta^\pm}$ are linear polynomials on $\Delta^\pm$, respectively, and satisfy
$$
\int_e v^+ \,ds =\int_e v^-\, ds.
$$
For boundary edges $e$, the edge integrals of $v\in V_\mathcal{T}$ vanish.
The discrete $H^1_0$ space associated with $\mathcal{T}$, 
containing both $V_{\mathcal{T}}$ and $H_0^1(\Omega)$ as closed subspaces and appropriate for the convergence analysis of the nonconforming P1 element method  \cite{BS}, consists of functions $u$ for which $u|_\Delta \in H^1(\Delta)$ for any triangle $\Delta\in \mathcal{T}$ and that satisfy the same compatibility conditions for integrals  along interior and boundary edges $e\in\mathcal{E}$ as indicated above for elements of $V_\mathcal{T}$. It will be denoted by $H^1_0(\Omega,\mathcal{T})$. The expression 
$$
|u|_{H^1,\mathcal{T}}:= \left(\sum_{\Delta\in\mathcal{T}} |u|_{1,\Delta}^2 \right)^{1/2},\qquad |u|_{1,\Delta}:=\left(\int_{\Delta} |\nabla u|^2\,dxdy\right)^{1/2},
$$
defines a norm on $H^1_0(\Omega,\mathcal{T})$ which turns it into a Hilbert space with scalar product denoted by $(\cdot,\cdot)_{H^1,\mathcal{T}}$. We call $
|u|_{H^1,\mathcal{T}}$ the discrete $H^1$ norm, and note that it coincides with the standard norm for $u\in H^1_0(\Omega)$. Here and throughout the paper, $\nabla u =(u_x,u_y)$ is the gradient of 
$u$, and $u_x,u_y,u_{xx},u_{x,y},u_{yy},\ldots$ is our notation for the partial derivatives of $u$ (if properly defined).
Thus, the variational problem of finding $u_{\mathcal{T}}\in V_{\mathcal{T}}$ such that
$$
(u_{\mathcal{T}},v)_{H^1,\mathcal{T}} = (f,v)_{L_2}\qquad \forall\; v\in V_{\mathcal{T}},
$$
has a unique solution which we call the Galerkin solution of (\ref{Poisson}) in $V_{\mathcal{T}}$.

In this paper we are concerned with estimates for the discrete $H^1$ norm error (called Galerkin error for short)
$$
E_{\mathcal{T}}(u):=|u-u_{\mathcal{T}}|_{H^1,\mathcal{T}} = \left(\sum_{\Delta\in\mathcal{T}} \int_\Delta |\nabla (u-u_{\mathcal{T}})|^2\,dxdy\right)^{1/2},
$$
if the solution $u$ of (\ref{Poisson}) is in $H^2(\Omega)$.
The second Strang Lemma implies that 
\begin{equation}\label{EBAC}
\max(E_{BA,\mathcal{T}}(u),E_{C,\mathcal{T}}(u))\le E_{\mathcal{T}}(u)\le E_{BA,\mathcal{T}}(u)+E_{C,\mathcal{T}}(u),
\end{equation}
i.e., that estimating the Galerkin error requires estimating both the best approximation error
$$
E_{BA,\mathcal{T}}(u) := \inf_{v\in V_{\mathcal{T}}}  |u-v|_{H^1,\mathcal{T}}
$$
of the solution $u$ by elements of $V_{\mathcal{T}}$, and the consistency error
$$
E_{C,\mathcal{T}}(u) := \sup_{w\in V_{\mathcal{T}}:\,|w|_{H^1,\mathcal{T}}=1} |(u,w)_{H^1,\mathcal{T}} -(f,w)_{L_2}|.
$$
In contrast to conforming P1 elements \cite{BA,Ku,Os}, in the nonconforming P1 case the best approximation error $E_{BA,\mathcal{T}}$ admits an
optimal bound for any $\mathcal{T}$. To formulate it, consider the Crouzeix-Raviart interpolation operator $P_{\mathcal{T}}:\;H^1_0(\Omega,\mathcal{T})\to 
V_{\mathcal{T}}$ introduced in \cite{CR} and defined by the condition
$$
\int_e (u-P_{\mathcal{T}}u)\,ds =0\qquad \forall\;e\in\mathcal{E}.
$$
The following result is a consequence of, e.g., Lemma 2.2 in \cite{AD}.
\begin{theo}\label{lem1} If $u\in H^2(\Omega)\cap H^1_0(\Omega)$ then, with a constant $C_0$ independent of $\mathcal{T}$ we have
$$
E_{BA,\mathcal{T}}(u) \le |u-P_{\mathcal{T}}u|_{H^1,\mathcal{T}} \le C_0 \left(  \sum_{\Delta\in\mathcal{T}} h_\Delta^2
|u|_{2,\Delta}^2        \right)^{1/2}\le C_0 h_{\mathcal{T}}|u|_{H^2},
$$
where 
$$
|u|_{2,\Delta}:=\left(\int_\Delta |D^2u|^2\,dxdy\right)^{1/2},\qquad |D^2u|^2:=u_{xx}^2+2u_{xy}^2+u_{yy}^2,
$$
and $|u|_{H^2}:=|u|_{2,\Omega}$ stands for the $H^2$ semi-norm of $u$.
\end{theo}

Unfortunately, the consistency error $E_{C,\mathcal{T}}(u)$ does not admit a similar estimate with constants uniform in $\mathcal{T}$.
Indeed, the standard estimate of $E_{C,\mathcal{T}}(u)$ is based on the transformation
\begin{equation}\label{Cerror}
(u,w)_{H^1,\mathcal{T}} -(f,w)_{L_2}= \sum_{e\in \mathcal{E}} \int_e (\nabla u\cdot n_{e}) [w]\,ds,\qquad w\in V_{\mathcal{T}},
\end{equation}
where $n_{e}$ is a fixed unit normal with respect to the edge $e$, and $[w]$ denotes the (properly signed) difference of the traces of $w$ from both sides
of $e$ (set $w=0$ outside $\Omega$). When each of these edge integrals is bounded by the trace theorem, see \cite{Br,BS}, a dependence on the shape of the triangles
attached to $e$ enters the constants. Implicitly, this can be seen from \cite[Theorem 6.2]{CGR} which contains the following estimate for the Galerkin error
(for simplicity, we do not state it with the explicit constants given in \cite{CGR}):
\begin{theo}\label{theo1} If $u\in H^2(\Omega)\cap H^1_0(\Omega)$ then, with constants $C_1$, $C_2$ independent of $\mathcal{T}$, we have
$$
E_{\mathcal{T}}(u) \le   \left( \sum_{\Delta\in\mathcal{T}} h_\Delta^2 \left\{C_1^2 \int_\Delta |f-\bar{f}_\Delta|^2\,dxdy +
C_2^2\tan^2(\frac{\alpha_{\Delta}}2) \int_\Delta |D^2u|^2\,dxdy\right\}       \right)^{1/2}
$$
$$
 \le h_{\mathcal{T}}(C_1\|f\|_{L_2} + C_2\tan(\frac{\alpha_{\mathcal{T}}}2) |u|_{H^2}),\qquad\qquad\qquad\qquad
$$
where $\bar{f}_\Delta:=|\Delta|^{-1}\int_\Delta f \,dxdy$ denotes the average value of $f$ on $\Delta$.
\end{theo}

The appearance of the factor $\tan(\alpha_{\mathcal{T}}/2)$ is troublesome, as it indicates a deterioration
of the error bound  if $\alpha_{\mathcal{T}}\to \pi$. Moreover, for sequences of triangulations
with $h_{\mathcal{T}}\tan(\alpha_{\mathcal{T}}/2) \to \infty$ even boundedness of the Galerkin error is not guaranteed! Whether $E_{\mathcal{T}}(u)\to \infty$
may happen for some $u\in H^2(\Omega)\cap H^1_0(\Omega)$ is doubtful but currently not disproved. This question is closely related to a possible deterioration of the constant in the discrete Friedrichs inequality
\begin{equation}\label{Fried0}
\|w\|_{L_2}\le C_{\Omega,\mathcal{T}}|w|_{H^1,\mathcal{T}}\qquad \forall\;w\in V_{\mathcal{T}},
\end{equation}
namely, if, for fixed polygonal $\Omega$, the supremum of the optimal constants $C_{\Omega,\mathcal{T}}$ in (\ref{Fried0}) over all possible
$\mathcal{T}$ may become infinity. There is some ambiguity on the dependence of $C_{\Omega,\mathcal{T}}$ on the shape regularity properties 
of $\mathcal{T}$ in the literature, see e.g. \cite{BS,Vo}, which we could not yet sort out. 

The family of triangulations $\mathcal{T}_{n,m}$ of the unit square, we concentrate on in this paper, does not exhibit such an extreme divergence behavior. However,
it shows that the dependency on $\alpha_{\mathcal{T}}$ present in the estimate of Theorem \ref{theo1} is essential, and that (bounded) divergence
of the nonconforming P1 method is possible. Let us introduce the notation 
used in Section \ref{sec3}. We consider the solution $u(x,y):= x(1-x)y(1-y)$ of the Poisson problem
\begin{equation}\label{Poisson1}
-\Delta u(x,y) = f(x,y):=2(x(1-x)+y(1-y)),\quad (x,y)\in [0,1]^2, 
\end{equation}
equipped with homogeneous Dirichlet boundary conditions
$$
u(0,y)=u(1,y)=u(x,0)=u(x,1)=0,\quad x,y\in [0,1], 
$$
and the associated sequence of nonconforming P1 element Galerkin solutions
$$
u_{n,m}:=u_{\mathcal{T}_{n,m}}\in V_{n,m}:=V_{\mathcal{T}_{n,m}}, \qquad m\ge n\ge 1.
$$
Even though Figure \ref{fig1} is self-explaining,
we give the formal definition of the triangulation $\mathcal{T}_{n,m}$. It is generated by the intersection of three line systems with $[0,1]^2$, namely
\begin{eqnarray*}
&&\{(x,y):\;y=\frac{j}{2m},\;x\in [0,1]\}_{j=1,\ldots,2m-1}, \\
&&\{(x,y):\;y=\frac{n}{m}x+ \frac{j}{m},\;x\in [0,1]\}_{j=1-m,\ldots,m-1},\\
&&\{(x,y):\;y=-\frac{n}{m}x+ \frac{j}{m},\;x\in [0,1]\}_{j=1,\ldots,2m-1}. 
\end{eqnarray*}
Its vertex set 
consists of all points $P_{i,j}=(\frac{i}{2n},\frac{j}{2m})$ with indices $i=0,2,\ldots,2n$ if $j=0,2,\ldots,2m$ is even, and indices $i=0,1,3,\ldots,2n-1,2n$ if $j=1,3,\ldots,2m-1$ is odd. The typical triangle
$\Delta$ in $\mathcal{T}_{n,m}$ has its longest edge of length $1/n$ located parallel to the $x$-axis, an associated height of length $1/(2m)$, area $|\Delta|=1/(4nm)$, and two remaining sides of equal length. It becomes severely distorted, with the maximum angle $\alpha_{\Delta}$ satisfying
$\tan(\alpha_{\Delta}/2)=m/n$, if $m/n\to \infty$ (the exceptional triangles along the vertical sides of the square are right-angled, have shorter longest edges,
and area $1/(8nm)$). Thus, we have
$$
h_{\mathcal{T}_{n,m}}=\frac1n,\qquad \tan(\frac{\alpha_{\mathcal{T}_{n,m}}}{2})=\frac{m}n, \qquad m\ge n\ge 1.
$$
The triangulations $\mathcal{T}_{n,m}$ have been used
in \cite{BA,Os} for studying $H^1$ best approximation with conforming P1 elements but seem to have appeared for the first time in H. Schwarz' seminal
note \cite{Sc} on the definition of the surface area by triangular approximation. 

We denote by $E_{n,m}=|u-u_{n,m}|_{H^1,\mathcal{T}_{n,m}}$ the Galerkin error of our model problem with respect to $\mathcal{T}_{n,m}$. Then Theorem \ref{theo1}
gives the upper bound
\begin{equation}\label{Enm0}
E_{n,m} \le C_3 \frac{m}{n^2},\qquad m\ge n\ge 1,
\end{equation}
where the constant $C_3$ is independent of $n$ and $m$. The main result of this paper is a two-sided estimate for $E_{n,m}$ and shows that the upper estimate (\ref{Enm0}) is essentially sharp in the range $n\le m\le n^2$.

\begin{theo}\label{theo2} For the model problem (\ref{Poisson1}) with solution $u(x,y)=x(1-x)y(1-y)$ and the family of triangulations $\mathcal{T}_{n,m}$ we have
\begin{equation}\label{Enm}
C'_4 \min(1,\frac{m}{n^2})\le E_{n,m} \le C_4 \min(1,\frac{m}{n^2}),\qquad m\ge n\ge 1,
\end{equation}
with constants $C_4,C'_4$ independent of $n$ and $m$. In particular, to achieve convergence in the discrete $H^1$ semi-norm for a certain sequence
of triangulations $\mathcal{T}_{n,m}$ with $n\to \infty$, one needs to satisfy $m/n^2\to 0$.
\end{theo}

The behavior of the Galerkin error for our model problem needs to be contrasted with the behavior of the best approximation error:
\begin{equation}\label{LowerBA}
E_{BA,\mathcal{T}_{n,m}}(u)\approx \frac1n,\qquad m\ge n\ge 1.
\end{equation}
The upper estimate in (\ref{LowerBA}) follows from Theorem \ref{lem1}, a matching lower bound is obtained if we invoke the two-sided Poincar\'{e} inequality
\begin{equation}\label{Poinc}
\inf_{c\in \mathbb{R}} \|v-c\|^2_{L_2(\Delta)} = \|v-\bar{v}_\Delta\|^2_{L_2(\Delta)}\approx \int_\Delta (\frac1{n^2}v_x^2+\frac1{m^2}v_y^2)\, dxdy, \qquad v\in H^1(\Delta),
\end{equation}
for the best approximation by constants, valid for any triangle $\Delta\in \mathcal{T}_{n,m}$ and any fixed polynomial $u(x,y)$ with positive constants depending
on the degree. To see (\ref{Poinc}), just use the coordinate transform $x'=x$,
$y'=\frac{m}{n}y$, apply the equivalence of $H^1$ semi-norm and $L_2$ norm on the finite-dimensional subspace of $H^1(\Delta')$ consisting of polynomials of fixed  degree with zero average
which holds, with uniform constants, for the transformed, undistorted triangle $\Delta'$, and then transform back. If one applies (\ref{Poinc}) separately to
the partial derivatives $u_x$ and $u_y$ of the solution of (\ref{Poisson1}), and adds the results for all $\Delta\in \mathcal{T}_{n,m}$,
then
\begin{eqnarray*}
E_{BA,\mathcal{T}_{n,m}}(u)^2&\ge& \sum_{\Delta\in\mathcal{T}_{n,m}} \inf_{c,c'\in \mathbb{R}} (\|u_x-c\|^2_{L_2(\Delta)}+
\|u_y-c'\|^2_{L_2(\Delta)})\\ 
&\ge& \frac{C'_0}{n^2}\int_{\Omega} (u_{xx}^2+u_{xy}^2+\frac{n^2}{m^2}u_{yy}^2)\,dxdy
\end{eqnarray*}
with some $C'_0>0$. This shows the lower bound in (\ref{LowerBA}). Thus, our main result formulated in Theorem \ref{theo2} is equivalent to showing a two-sided estimate similar to (\ref{Enm}) for the consistency error $E_{C,\mathcal{T}_{n,m}}(u)$.

\section{Proof of Theorem \ref{theo2}}\label{sec3}

We first deal with the upper bound in (\ref{Enm}). Due to (\ref{Enm0}) all we need is to establish a complementing upper bound for $E_{n,m}$ by a constant, independent of $n$ and $m$. Since
$$
E_{n,m}\le E_{BA,\mathcal{T}_{n,m}}(u)+E_{C,\mathcal{T}_{n,m}}(u)\le \|u\|_{H^1} +E_{C,\mathcal{T}_{n,m}}(u),
$$
and
$$
|(u,w)_{H^1,\mathcal{T}_{n,m}} -(f,w)_{L_2}|\le \|u\|_{H^1}|w|_{H^1,\mathcal{T}_{n,m}}+\|f\|_{L_2}\|w\|_{L_2},\qquad w\in V_{n,m},
$$
the upper bound in (\ref{Enm}) holds with constant $C_4=\max(C_3,2\|u\|_{H^1}+\frac12 \|f\|_{L_2})$, since for the triangulations $\mathcal{T}_{n,m}$
we have the discrete Friedrichs inequality
\begin{equation}\label{Fried}
\|w\|_{L_2} \le \frac12 |w|_{H^1,\mathcal{T}_{n,m}},\qquad w\in V_{n,m},\qquad m\ge n\ge 1.
\end{equation}
Since we could not find a reference for (\ref{Fried}) in the literature, we give the elementary argument in Section \ref{sec51}. 

The rest of the proof is concerned with proving the matching lower bound in (\ref{Enm}). As was pointed out before, this is equivalent to establishing the appropriate lower bound for 
$$
E_{C,\mathcal{T}_{n,m}}(u) = \sup_{0\neq w\in V_{n,m}} \frac{|(u,w)_{H^1,\mathcal{T}_{n,m}} -(f,w)_{L_2}|}{|w|_{H^1,\mathcal{T}_{m,n}}}.
$$
To this end, it is enough to pick a suitable $\tilde{w}\in V_{n,m}$, estimate its discrete $H^1$ norm from above, the 
consistency term $|(u, \tilde{w})_{H^1,\mathcal{T}_{n,m}} -(f,\tilde{w})_{L_2}|$ from below, and check the quotient of these estimates. We arrived at a good guess for a such a candidate $\tilde{w}$ after performing some numerical experiments, see Section \ref{sec4}. 
We define the nodal values $\tilde{w}(M_e)$ as follows: For all edges $e$ in the lower left subsquare $\Omega':=[0,\frac12]^2$ of $\Omega$, we set 
\begin{equation}\label{DefW}
\tilde{w}(M_e)=\left\{\begin{array}{ll} 0, &  \mbox{$e$ on the boundary, or parallel to the $x$-axis},\\
\psi(M_e),& \mbox{$e$ has slope $n/m$},\\  -\psi(M_e),& \mbox{$e$ has slope $-n/m$}, \end{array}\right.
\end{equation}
where $\psi(x,y)=2hx u_{xy}(x,y)=2x(1-2x)(1-2y)$. Nodal values for the remaining part of $\Omega$ are obtained by symmetry, i.e., such that $\tilde{w}(1-x,y)=\tilde{w}(x,1-y)=
\tilde{w}(1-x,1-y)=\tilde{w}(x,y)$ for all $(x,y)\in \Omega'$. Note that this $\tilde{w}$ is highly oscillating, and related to the mixed derivative
$u_{xy}=(1-2x)(1-2y)$, scaled by $h$ and with values damped towards the vertical edges of $\Omega$ by the factor $\min(2x,2(1-x))$.

\begin{figure}[htb] 
\begin{center}
\includegraphics[width=1.2\textwidth]{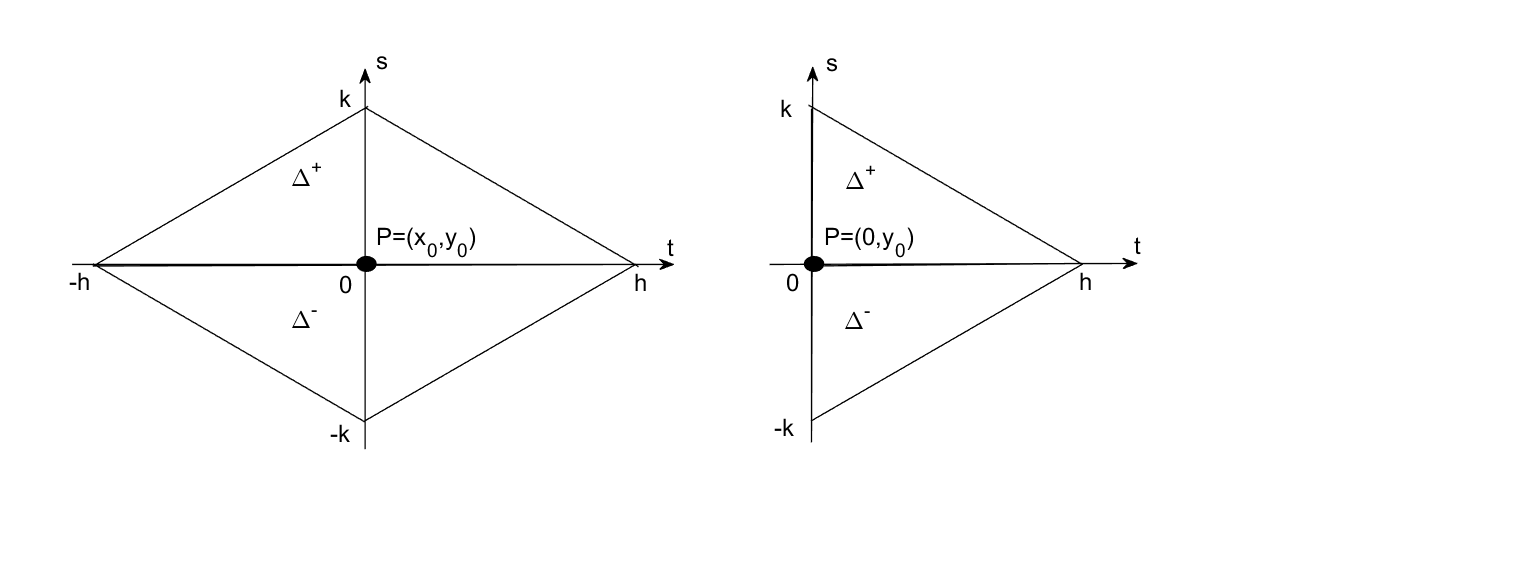}
\caption{Triangle pairs in the interior (on the left), and attached to the boundary (on the right)}
\label{fig2}
\end{center}
\end{figure}

By symmetry, we need to evaluate the integrals defining $(u,\tilde{w})_{H^1,\mathcal{T}_{n,m}}$, $(f,\tilde{w})_{L_2}$, and $\|\tilde{w}\|^2_{H^1,\mathcal{T}_{n,m}}$ only on the subsquare $\Omega'$. Thus,  estimates will be conducted for the triangles depicted in Figure \ref{fig2} that intersect with $\Omega'$. We use the notation introduced by Figure \ref{fig2}, with the reference point $P=(x_0,y_0)$ (resp. $P=(0,y_0)$) representing the origin of a local coordinate system $(t,s)$, and $h:=1/(2n)$ and $k:=1/(2m)$ the typical lengths in $t$- and $s$-direction, respectively. We also denote 
$$
\kappa:=h^2k^{-1}=\frac{m}{2n^2}. 
$$

Formulas for the piecewise constant gradient $\nabla \tilde{w}|_{\Delta^\pm} =: (\tilde{w}^\pm_x,\tilde{w}^\pm_y)$ for all triangles intersecting with $\Omega'$
follow from the definition of $\tilde{w}$ by elementary calculus, and immediately lead to estimates for the discrete
$H^1$ norm of $\tilde{w}$. The result is collected into the following lemma, see Section \ref{sec52} for its derivation.
\begin{lem}\label{lem2} Let $\tilde{w}\in V_{n,m}$ be given by (\ref{DefW}).\\
a) For the triangles $\Delta^{\pm}\subset \Omega'$ with reference point $P=(0,y_0)$, $0<y_0<1/2$ (see Figure \ref{fig2} on the right), we have
\begin{equation}\label{WxyB}
\tilde{w}_x^{\pm}=2h(1-h)(\pm(1-2y_0)+k),\qquad \tilde{w}^\pm_y=-2\kappa (1-h)((1-2y_0)\mp k),
\end{equation}
For the triangles $\Delta^{\pm}\subset \Omega'$ with reference point $P=(x_0,y_0)$, $0<x_0<1/2$, $0\le y_0 \le 1/2$ (see Figure \ref{fig2} on the left), we have
\begin{equation}\label{Wxy}
\tilde{w}_x^\pm=(4x_0(1-2x_0)-2h^2)(\mp(1-2y_0)+k),\qquad \tilde{w}_y^\pm=\kappa(4x_0-1)(1-2y_0\mp k).
\end{equation}
Finally, for the triangles $\Delta^{\pm}$ with reference point $P=(1/2,y_0)$, $0\le y_0\le 1/2$, on the symmetry line $x=1/2$, we have
\begin{equation}\label{WxyM}
\tilde{w}_x^\pm=0,\qquad \tilde{w}_y^{\pm} =2\kappa(1-h)(1-2y_0\mp k).
\end{equation}
b) The discrete $H^1$ norm of $\tilde{w}$ satisfies
\begin{equation}\label{WBound}
|\tilde{w}|_{H^1,\mathcal{T}_{n,m}}= \mathrm{O}(1+\kappa),\qquad m\ge n\ge 1.
\end{equation}
\end{lem}

We come to the lower estimate for the consistency term evaluated at $\tilde{w}$. As it turns out, the dominating contributions to the consistency term
come from the integrals
$$
\int_{\Delta^\pm} u_y\tilde{w}_y\, dxdy
$$
for interior triangle pairs $\Delta^{\pm}\subset \Omega'$, as depicted in Figure \ref{fig2} on the left, and are of the order $m/n^2$. Other terms are negligible compared to them.
In particular, we have the following lemma whose proof is given in Section \ref{sec53}.
\begin{lem}\label{lem3} For the $\tilde{w}$ under consideration and the right-hand side $f$ in (\ref{Poisson1}), we have
\begin{equation}\label{fw}
|(f,\tilde{w})_{L_2}| =\mathrm{O}(kh^2), \qquad m\ge n\ge 1.
\end{equation}
\end{lem}

The crucial part of the proof is a lower bound for $(u,\tilde{w})_{H^1,\mathcal{T}_{n,m}}$. We first deal with the contributions to $(u,w)_{H^1,\mathcal{T}_{n,m}}$ from the triangles $\Delta^\pm\subset \Omega'$  depicted in Figure \ref{fig2} on the left. Have in mind that
in local coordinates we have 
$$
u_x(x_0+t,y_0+s)=(1-2x_0-2t)(y_0(1-y_0)+(1-2y_0)s-s^2),
$$ analogously for $u_y(x_0+t,y_0+s)$, while $\tilde{w}_x^\pm$, $\tilde{w}_y^\pm$ are constant
on $\Delta^\pm$, respectively, and given by (\ref{Wxy}). Using the simplifications based on symmetry arguments and integration over triangles as detailed in Section \ref{sec53}, we have
\begin{eqnarray*}
\int_{\Delta^{\pm}} u_x \,dxdy &=& (1-2x_0) \int_{\Delta^{\pm}_0} y_0(1-y_0)+(1-2y_0)s-s^2 \,dtds\\
&=&hk(1-2x_0)\left(y_0(1-y_0)\pm \frac13 (1-2y_0)k-\frac16 k^2\right)\\
&=&hk(1-2x_0)(y_0(1-y_0) + \mathrm{O}(k)),
\end{eqnarray*}
and 
\begin{eqnarray*}
\int_{\Delta^{\pm}} u_y \,dxdy &=& \int_{\Delta^{\pm}_0} (x_0(1-x_0)-t^2)(1-2y_0-2s)\,dtds\\
&=&hk((x_0(1-x_0)-\frac16 h^2)(1-2y_0)\mp\frac23 kx_0(1-x_0) \pm \frac{1}{15}h^2k)\\
&=&hk((1-2y_0)x_0(1-x_0) +\mathrm{O}(k+h^2)).
\end{eqnarray*}
Here $\Delta^{\pm}_0$ denotes the triangle pair associated with reference point $(0,0)$. Substituting the values 
$$
\tilde{w}^\pm_x=\mp 4x_0(1-2x_0)(1-2y_0)+\mathrm{O}(h^2+k),\qquad \tilde{w}^\pm_y = \kappa((4x_0-1)(1-2y_0)+\mathrm{O}(k)),
$$ 
obtained from (\ref{Wxy}), we get
\begin{eqnarray*}
\int_{\Delta^{\pm}} \nabla u\cdot \nabla \tilde{w}\,dxdy &=& \tilde{w}^\pm_x \int_{\Delta^{\pm}} u_x \,dxdy +\tilde{w}^\pm_y \int_{\Delta^{\pm}} u_y \,dxdy\\
&=& \mp 4 hk (x_0(1-2x_0)^2y_0(1-y_0)(1-2y_0) + \mathrm{O}(h^2+k))\\
&& \qquad + hk\kappa((4x_0-1)x_0(1-x_0)(1-2y_0)^2  + \mathrm{O}(h^2+k)).
\end{eqnarray*}
If we sum with respect to the $\mathrm{O}(nm)$ triangles in $\Omega'$  considered so far (call the result $\Sigma'$), we see that
\begin{equation}\label{uw1}
\Sigma' = \kappa (I'+\mathrm{O}(h+k^2/h^2)),
\end{equation}
with a constant $I'>0$ given below.
Indeed, for the terms in the sum $\Sigma'$ related to the gradient in $x$-direction, the leading parts $\mp 4 hk x_0(1-2x_0)^2y_0(1-y_0)(1-2y_0)$ cancel for triangle pairs $\Delta^\pm\subset \Omega'$ with the same  reference point $P=(x_0,y_0)$, and vanish for triangles $\Delta^+$ with $y_0=0$ and $\Delta^-$ with $y_0=1/2$,
respectively. Therefore, only the subdominant part $\mathrm{O}(hk(h^2+k))$ needs to be taken into account which gives an overall $\mathrm{O}(h^2+k)=
\mathrm{O}(\kappa(k+k^2/h^2))$
contribution to $\Sigma'$. Moreover, for the terms in $\Sigma'$ related to the gradient in $y$-direction, the sum of the leading factors 
$hk(4x_0-1)x_0(1-x_0)(1-2y_0)^2$ (without the factor $\kappa$) tends to the integral 
$$
I':=\int_{\Omega'} (4x-1)x(1-x)(1-2y)^2 \,dxdy  = \frac1{384}
$$
at a speed of at least $\mathrm{O}(h)$ as $h,k\to 0$. Altogether, this gives (\ref{uw1}) if one takes the common factor $\kappa=h^2/k$ out, and uses $k=\mathrm{O}(h)$. 
We can silently include into $\Sigma'$ the contributions from the $\mathrm{O}(m)$ triangles $\Delta^\pm$ crossing the symmetry line $x=1/2$, as the 
estimation steps are identical, with the only change that (\ref{Wxy}) is replaced by (\ref{WxyM}).

The contribution of the remaining triangles $\Delta^\pm$ with $P=(0,y_0)$, depicted in Figure \ref{fig2} on the right and attached to the left boundary of $\Omega'$, is negligible compared to the leading part in the lower estimate (\ref{uw1}). Indeed, we again expand in local coordinates $(t,s)$ as
$$
u_x(t,y_0+s)=-2t(y_0(1-y_0)+(1-2y_0)s-s^2),\qquad u_y(t,y_0+s) = (t-t^2)(1-2y_0-2s),  
$$
where $0\le t\le h(1-|s|/k)$, $0\le s\le k$ for $\Delta^+$, and $-k\le s\le 0$ for $\Delta^-$, respectively, and compute with (\ref{BetaL}) the integrals
$$ 
\int_{\Delta^\pm} u_x\,dxdy =-\frac{h^2k}3\left(y_0(1-y_0)\pm  \frac{k}2(1-2y_0)-\frac{k^2}5\right)=\mathrm{O}(h^2k),
$$
and, similarly, 
$$
\int_{\Delta^\pm} u_y\,dxdy =\frac{h^2k}6 \left((1-2y_0)(1-\frac{h}2)\mp k( \frac{1}2 - \frac{h}5)\right)=\mathrm{O}(h^2k).
$$
Combining this with 
$$
\tilde{w}_x^\pm=\pm 2h(1-2y_0+\mathrm{O}(h)),\qquad \tilde{w}_y^\pm= -2\kappa (1-2y_0+\mathrm{O}(h)),
$$
see (\ref{WxyB}), we obtain the rough estimates
$$
\tilde{w}_x^\pm\int_{\Delta^\pm} u_x\,dxdy=\mathrm{O}(h^3k)=\mathrm{O}(\kappa hk^2),\qquad 
\tilde{w}_y^\pm\int_{\Delta^\pm} u_y\,dxdy=\mathrm{O}(\kappa h^2k).
$$
Summing the contributions with respect to all $\mathrm{O}(m)$ triangles attached to the boundary $x=0$ of $\Omega'$ (call the result $\Sigma''$),
we get 
\begin{equation}\label{uw2}
\Sigma'' = \mathrm{O}(\kappa h^2).
\end{equation}

Combining (\ref{uw1}), (\ref{uw2}), and (\ref{fw}), we see that 
\begin{equation}\label{lb}
(u,\tilde{w})_{H^1,\mathcal{T}_{n,m}}-(f,\tilde{w})_{L_2} = 4(\Sigma'+\Sigma'')-(f,\tilde{w})_{L_2}
= \kappa (I' + \mathrm{O}(h+k^2/h^2)). 
\end{equation}
Eventually, by (\ref{lb}) and (\ref{WBound}), we get, with an absolute constant $C'_5>0$,
$$
E_{C,\mathcal{T}_{n,m}}(u) \ge \frac{(u,\tilde{w})_{H^1,\mathcal{T}_{n,m}}-(f,\tilde{w})_{L_2}}{|\tilde{w}|_{H^1,\mathcal{T}_{n,m}}}
\ge C'_5\frac{\kappa}{1+\kappa}\ge\frac{C'_5}3 \min(1,m/n^2),
$$
 if $m\ge n\ge n_0$ with $n_0$ large enough, and $n/m\le \epsilon_0$ with $\epsilon_0$ small enough. This proves (\ref{Enm}) in the asymptotic range.
For the remaining values $m\ge n$, note that for them $\tan(\alpha_{\mathcal{T}_{n,m}}/2)=m/n\le C_6$ for some absolute $C_6$ depending on 
$n_0$, and $\epsilon_0$, i.e., these remaining triangulations $\mathcal{T}_{n,m}$ uniformly satisfy the maximum angle condition. Thus, in this case 
$1/n\le m/(C_6n^2)$, and the lower bound in (\ref{Enm}) is taken care of by the lower bound (\ref{LowerBA}) for the best discrete $H^1$ approximation error of our $u$. With the constant $C'_4$ in (\ref{Enm}) defined from $C'_5$, $C_6$, and from the constant in (\ref{LowerBA}) in a proper way, Theorem \ref{theo2} is now fully proved.

\section{Numerical Examples and Further Remarks}\label{sec4}
\begin{figure}
\centering
\includegraphics[width=0.45\textwidth,height=0.43\textwidth]{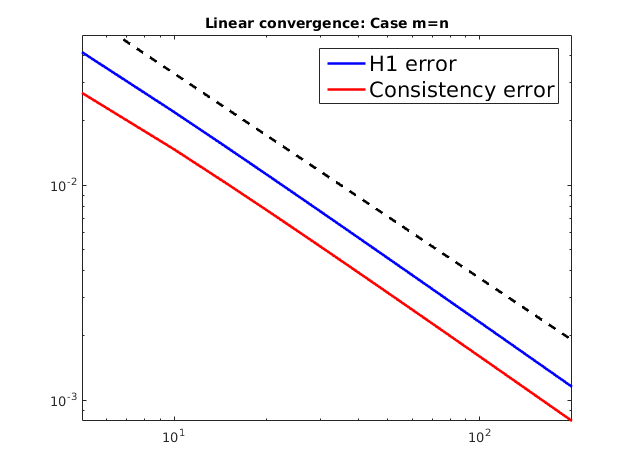}
\includegraphics[width=0.45\textwidth,height=0.43\textwidth]{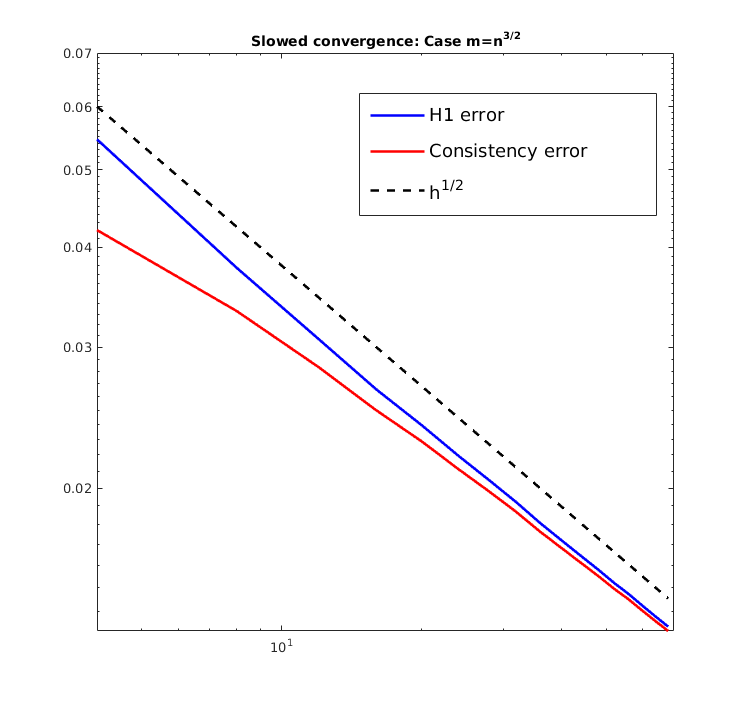}
\caption{Behavior of $E_{n,m}$ for $m=n$ (optimal order convergence, on the left) and for $m\approx n^{3/2}$ (slowed convergence, on the right)}
\label{fig3}
\end{figure}

\begin{figure}
\centering
\includegraphics[width=0.45\textwidth,height=0.43\textwidth]{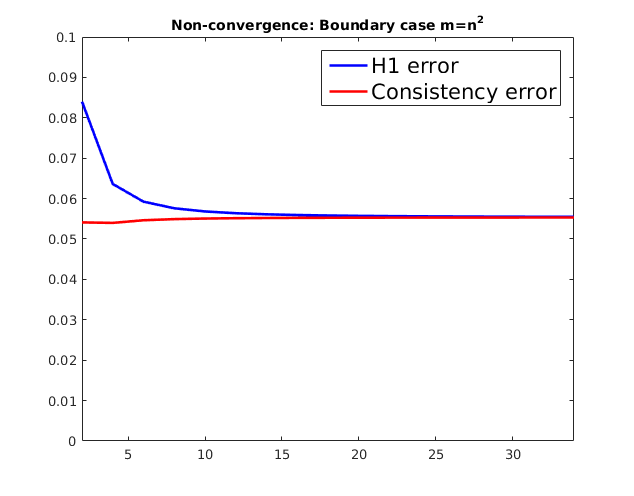}
\includegraphics[width=0.45\textwidth,height=0.43\textwidth]{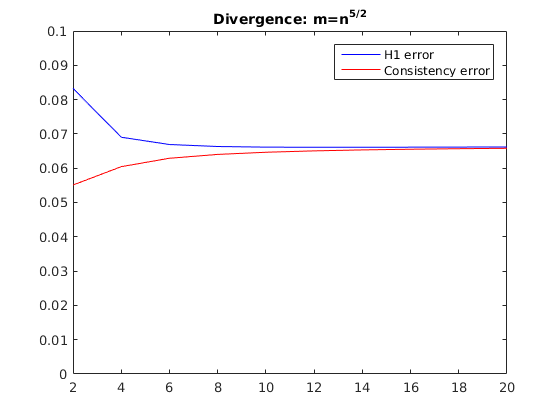}
\caption{Failure of convergence for $m=n^2$ (on the left) and for $m\approx n^{5/2}$ (on the right)}
\label{fig4}
\end{figure}

We have conducted a couple of numerical experiments in the pre-asymptotic range (with relatively small values $n$, $m$), for exactly the model problem 
described in the previous sections. We have used the standard nodal basis $\{\phi_e\}$ for nonconforming P1 elements associated with the interior edges of
$\mathcal{T}_{n,m}$, and computed the integrals defining the entries of the stiffness matrix $A$ and load vector $b$, as well as the 
error measures exactly (within machine accuracy).
First we  confirmed the result of Theorem \ref{theo2} by running simulations for values $m=n$,
$m\approx n^{3/2}$, $m=n^2$, and $m\approx n^{5/2}$, respectively, for a suitable range of values $n$. The first two cases shown in Figure \ref{fig3} illustrate optimal O$(n^{-1})$
and slowed O$(n^{-1/2})$ convergence, in agreement with (\ref{Enm}). The latter two cases demonstrate the failure of convergence if 
$m/n^2$ does not converge to zero, see Figure \ref{fig4}. Blue lines represent the Galerkin error, red lines the consistency error.

\begin{figure}
\centering
\includegraphics[width=0.45\textwidth,height=0.43\textwidth]{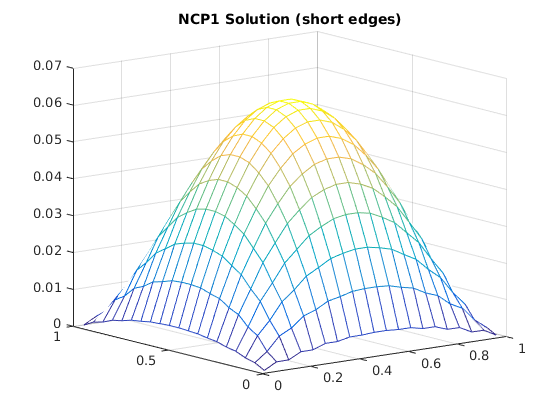}
\includegraphics[width=0.45\textwidth,height=0.43\textwidth]{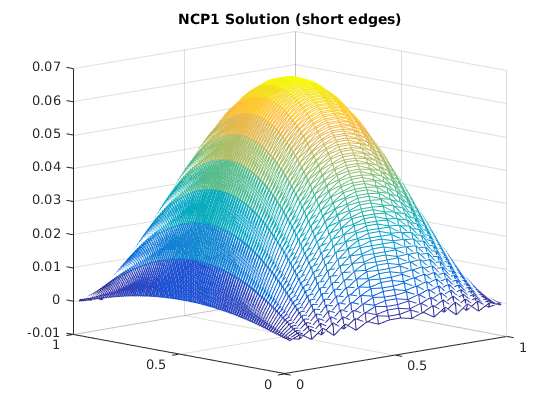}\\
\includegraphics[width=0.45\textwidth,height=0.43\textwidth]{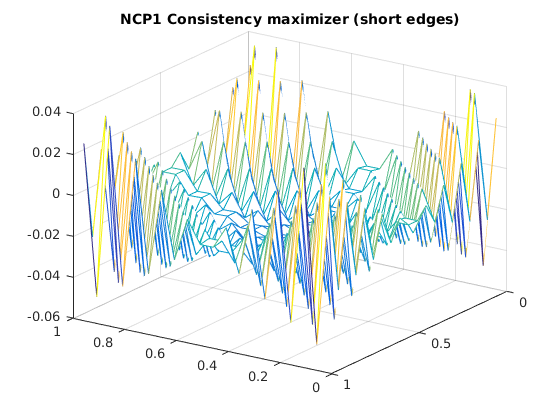}
\includegraphics[width=0.45\textwidth,height=0.43\textwidth]{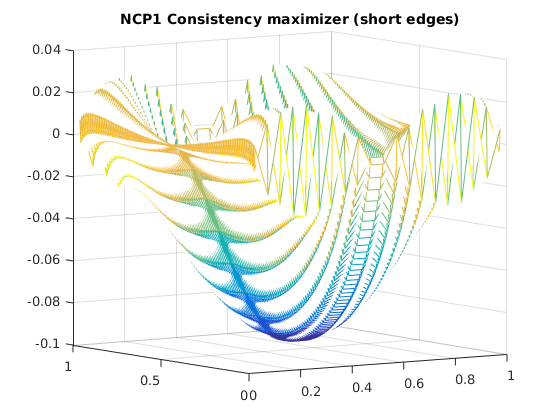}
\caption{Solutions (upper row) and maximizers for the consistency error (lower row) for $n=10$ and $m=n$ (on the left)
and $m=n^2$ (on the right)}
\label{fig5}
\end{figure}

We also needed some intuition on how an appropriate candidate $\tilde{w}$ for maximizing the consistency error should look like.
Since the constrained problem
$$
(u,\tilde{w})_{H^1,\mathcal{T}_{n,m}}-(f,\tilde{w})_{L_2}\to \max\quad \mbox{ subject to }\,|\tilde{w}|_{H^1,\mathcal{T}_{n,m}}=1
$$
is easy to solve, the coefficient vector of the maximizer $\tilde{w}$ and the value of $E_{C,\mathcal{T}_{n,m}}(u)$ can be found  from the formulas
$$
\tilde{x}=\pm A^{-1}(c-b)/\sqrt{(c-b)^TA^{-1}(c-b)},\qquad E_{C,\mathcal{T}_{n,m}}(u) = \sqrt{(c-b)^TA^{-1}(c-b)},
$$
where $c$ has entries $c_e=(u,\phi_e)_{H^1,\mathcal{T}_{n,m}}$. The result is visualized in Figure \ref{fig5} by depicting the nodal values of the Galerkin solution $u_{n,m}$
given by $x=A^{-1}b$  (upper row), and of the maximizer of the consistency error given by $\tilde{x}$ (lower row)
at the midpoints of edges with slope $\pm n/m$. We show two cases: $n=m=10$ (on the left), and  $n=10$, $m=n^2=100$ (on the right).
The graphs suggested a distinct oscillation behavior for $\tilde{w}$ which we slightly simplified to the choice for $\tilde{w}$ used in the
proofs of the previous section (it took us a while to realize that for the deterioration of the consistency error the non-oscillating 
part of $\tilde{w}$ visible in Figure \ref{fig5} is not essential). It also looks as if $u_{n,m}$ is still close to $u$ in $L_2$ and $L_\infty$ distance,
even in cases when the discrete $H^1$ error does not converge to zero. This is in contrast to the counterexamples for conforming P1 elements used in
\cite{Os}. 

Our example automatically provides similar matching lower bounds for lowest-order Raviart-Thomas
elements \cite{RT} if the mixed formulation of (\ref{Poisson}) is used. Indeed, due to \cite{Ma}, on each triangle $\Delta\in \mathcal{T}_{n,m}$,
the discrete flux $\sigma_{n,m}$ of the mixed method
belonging to the lowest-order Raviart-Thomas space on $\mathcal{T}_{n,m}$ and the gradient of the nonconforming P1 Galerkin solution $\tilde{u}_{n,m}$
of a modified Poisson problem with solution $\tilde{u}$ and piecewise constant right-hand side $\tilde{f}$ defined by
$$
(\tilde{f}|_\Delta)(x,y) =\bar{f}_\Delta, \qquad (x,y)\in \Delta,\qquad \Delta\in \mathcal{T}_{n,m},
$$ 
are related by
\begin{equation}\label{Marini}
\nabla \tilde{u}_{n,m}(x,y)-\sigma_{n,m}(x,y) = \frac12 \bar{f}_\Delta ((x,y)-M_\Delta),\qquad (x,y)\in \Delta.
\end{equation}
Here $M_\Delta$ denotes the barycenter of $\Delta$, and $\bar{f}_\Delta$ is the average value of $f$ on $\Delta$ as defined before. See also \cite{Br1,CGR,CPS}, where the connections between energy norm errors for conforming and nonconforming P1 elements as well as lowest-order Raviart-Thomas elements have been examined in order to obtain sharp a posteriori estimates for the Poisson problem. For our model problem (\ref{Poisson1}), since
on each triangle $\delta$ of $\mathcal{T}_{n,m}$
$$
\nabla u_{n,m} - \sigma_{n,m} = (\nabla u_{n,m} -  \nabla \tilde{u}_{n,m}) + (\nabla\tilde{u}_{n,m}-\sigma_{n,m}),
$$
the error of the lowest order mixed Raviart-Thomas method, i.e., the $L_2$ vector norm of $\nabla u -\sigma_{n,m}$, and the Galerkin error of the nonconforming P1 method is bounded by the sum of two  terms, namely $\|\nabla \tilde{u}_{n,m}-\sigma_{n,m}\|_{L_2}$ and 
$\|\nabla \tilde{u}_{n,m}-\nabla u_{n,m}\|_{L_2}$. Using (\ref{Marini}), the first one can be estimated by
$$
\frac12\left( \sum_{\Delta\in \mathcal{T}_{n,m}}\bar{f}_\Delta^2 \int_\Delta |(x,y)-M_\Delta|^2\,dxdy\right)^{1/2}
=\mathrm{O}(h_{\mathcal{T}_{m,n}}\|f\|_{L_2})=\mathrm{O}( \frac{1}n|u|_{H^2}),
$$
while for the second one
$$
\|\nabla \tilde{u}_{n,m}-\nabla u_{n,m}\|_{L_2}\le \|\nabla (\tilde{u}-u)\|_{L_2}\le C_7\|\tilde{f} -f\|_{L_2}\le C'_7 \frac{1}n
$$
due to the orthogonality properties of the nonconforming P1 Galerkin projection, 
elliptic regularity, and the fact that our $f$ is smooth.
Thus, from (\ref{Enm}) we conclude that
\begin{equation}\label{RTBound}
C'_8\min(1,m/n^2)\le \|\nabla u -\sigma_{n,m}\|_{L_2}\le C_8\min(1,m/n^2)
\end{equation}
with some positive constants $C'_8$, $C_8$, where the lower bound is guaranteed to hold if $n/m$ is small enough, i.e., when the maximum angle condition fails. 

What we did not consider in this note are extensions along the lines of \cite{Ku} where it was observed that long chains of distorted triangles are the reason for convergence deterioration in the conforming P1 case. For higher-order elements, similar effects are to be expected, even though there are differences (e.g., the critical exponent $\beta$ for which $m/n^\beta\not\to 0$ implies convergence failure grows with the polynomial degree).

We conclude with a sketch of the argument for a statement made in the introduction of this paper. Consider the model Poisson problem 
$$
-\Delta u(x,y) =1, \qquad (x,y)\in [0,1]^2,
$$
satisfying the boundary conditions $u(0,y)=u(1,y)=0$ in $x$-direction and periodic boundary conditions in $y$-direction
whose solution is given by the univariate polynomial $u(x)=\frac12 x(1-x)$.  In \cite{Os}, we showed that for this problem the conforming P1 element Galerkin solutions $\hat{u}_{n,m}$ on the triangulations $\mathcal{T}_{n,m}$ satisfy 
$$
\| u-\hat{u}_{n,m}\|_{H^1} \approx \min(1,m/n^2),\qquad m\ge n     \ge 1.
$$
I.e., for conforming P1 elements the $H^1$ energy norm convergence rate may degenerate with the mesh distortion even for an essentially one-dimensional solution. Interestingly enough, for this problem the nonconforming P1 Galerkin solutions converge at optimal speed:
\begin{equation}\label{P1NC0}
|u-u_{n,m}|_{H^1,\mathcal{T}_{n,m}} = \mathrm{O}(\frac1n) ,\qquad m\ge n     \ge 1.
\end{equation}
This also shows that the convergence behavior of conforming and nonconforming P1 Galerkin solutions may be drastically different if the triangulations violate the maximum angle condition.

To prove (\ref{P1NC0}), it is sufficient to bound the consistency error. In this case, it is convenient to use (\ref{Cerror}) and we will give the estimate
for any $C^2$ smooth $u=u(x)$ depending only on $x$. Indeed, all integrals in (\ref{Cerror}) with respect to horizontal edges of the triangulations $\mathcal{T}_{n,m}$ automatically vanish, since $\nabla u = (u'(x),0)$ and in this case $n_e=(0,\pm 1)$. Integrals with respect to the vertical boundary edges on the line $x=0$ also vanish,
since $\nabla u \cdot n_e=\pm u'(0)$ is constant and $[w]$ has zero average on $e$,
similarly for the vertical boundary edges on the line $x=1$. On all remaining edges, we have 
$$
\nabla u \cdot n_e = \pm  \frac{k}{\sqrt{h^2+k^2}}u'(x),
$$
while
$$
[w] = \pm \frac{2(x-x_e)\sqrt{h^2+k^2}}{h} (w_e^+- w_e^-),
$$
where $x_e$ is the $x$-coordinate of the midpoint $M_e$ of the edge $e$, and $w_e^\pm$ denotes the constant derivatives in direction $e$ of the restrictions of $w\in V_{n,m}$ to the two triangles $\Delta_e^\pm$ attached to $e$, respectively. These formulas can be checked by elementary calculus, the signs in them depend on the ordering of triangles and the choice of edge normals but are irrelevant for the subsequent
estimates. What is important is that $[w]$ has average zero on $e$, thus
\begin{eqnarray*}
\left|\int_e (\nabla u\cdot n_e)  [w]\, ds\right| &\le& \inf_c \frac{2k}{h}\left|\int_e  (u'(x)-c)  (x-x_e)(w_e^+- w_e^-)\, ds\right|\\
&\le& \frac{2k}{h} \frac{h}2\|u''\|_{L_\infty} \frac{h|e|}2 (|w_e^+| + |w_e^-|) \\
&\le& Ckh^2 (|\nabla(w|_{\Delta_e^+})| + |\nabla(w_{\Delta_e^-})|)
\end{eqnarray*}
with a constant $C$ depending on $u(x)$ only. Now apply the Cauchy-Schwarz inequality to the sum of these upper estimates. This gives
\begin{eqnarray*}
|(u,w)_{H^1,\mathcal{T}_{n,m}} -(f,w)_{L_2}|^2&\le& C \left({\sum_e}' h^3 k\right) \left({\sum_e}' hk (|\nabla(w|_{\Delta_e^+})|^2 + |\nabla(w_{\Delta_e^-})|^2)\right)\\
&\le& \frac{C'}{n^2} |w|_{H^1,\mathcal{T}_{n,m}}^2,
\end{eqnarray*}
where $\sum'_e$ 
is the sum over the O$(nm)$ edges with nontrivial edge integrals in (\ref{Cerror}), and $C'$ is a new absolute constant. 
This is the desired bound for the consistency error, and together with Theorem \ref{lem1} implies (\ref{P1NC0}).

\bibliographystyle{amsplain}

\begin{thebibliography}{ABC}
\bibitem {AD} G. Acosta, R. G. Dur\'{a}n, {The maximum angle condition for mixed and non conforming elements. {A}pplication to the Stokes equations}, SIAM J. Numer. Anal.  \textbf{37} (1999), 18--36. Zbl 0948.65115,
\bibitem {BA} I. Babu\v{s}ka, A. K. Aziz, {On the angle condition in the finite element method}, SIAM J. Numer. Anal. \textbf{13} (1976), 214--226. Zbl 0324.65046,
\bibitem{Br} D. Braess, {Finite Elements. Theory, Fast Solvers and Applications in Solid Mechanics}, 3rd  ed., Cambridge Univ. Press, 2007.
Zbl 1118.65117, 
\bibitem{Br1} D. Braess, {An a posteriori error estimate and a comparison theorem for the nonconforming $P_1$ element}, Calcolo \textbf{46} (2009), 149--155.
Zbl 1192.65142,
\bibitem{SB} S. C. Brenner, Poincar\'{e}-Friedrichs inequalities for piecewise $H^1$ functions, SIAM J. Numer. Anal. \textbf{41} (2003), 306--324. Zbl 1045.65100,
\bibitem {BS} S. C. Brenner, L. R. Scott, {The Mathematical Theory of Finite Element Methods}, Texts in Appl. Math., vol. 15, Springer, New York, 2008. Zbl 0804.65101, 
\bibitem {CGR} C. Carstensen, J. Gedicke, D. Rim, {Explicit error estimates for Courant, Crouzeix-Raviart and Raviart-Thomas finite element methods},
J. Comput. Math. \textbf{30:4} (2012), 337--353. Zbl 1274.65290,
\bibitem {CPS} C. Carstensen, D. Peterseim, M. Schedensack, {Comparison results on finite element methods for the Poisson model problem}, SIAM J. Numer. Anal. \textbf{50} (2012), 2803--2823. Zbl 1261.65115,
\bibitem{CR} M. Crouzeix, P.-A. Raviart, {Conforming and nonconforming finite element  methods for solving the stationary Stokes equations} I, R.A.I.R.O. Anal. Num\'{e}r. \textbf{7:3} (1973), 33--76. Zbl 0302.65087,
\bibitem{HKK} A. Hannukainen, S. Korotov, M. K\v{r}{\'\i}\v{z}ek, {The maximum angle condition is not necessary for
convergence of the finite element method}, Numer. Math. \textbf{120} (2012), 79--88.
Zbl 1255.65196, 
\bibitem{Ja} P. Jamet, {Estimations d$'$erreur pour des \'{e}l\'{e}ments finis droits presque d\'{e}g\'{e}n\'{e}r\'{e}s}, R.A.I.R.O. Anal.
Num\'{e}r. \textbf{10} (1976), 43--61. Zbl 0346.65052,
\bibitem{Ku} V. Ku\v{c}era, {On necessary and sufficient conditions for finite element convergence}, arXiv:1601.02942 (2016). 
\bibitem{Ma} L. D. Marini, {An inexpensive method for the evaluation of the solution of the lowest order Raviart-Thomas mixed method}, SIAM J. Numer. Anal. \textbf{22:3} (1985), 493--496. Zbl 0573.65082,
\bibitem {Os} P. Oswald, {Divergence of FEM: Babu\v{s}ka-Aziz triangulatiuons revisited}, Appl. Math. \textbf{60:5} (2015), 473--484. Zbl 1363.65202,
\bibitem{RT} P.-A. Raviart, J. M. Thomas, {A mixed finite element  method for second order elliptic problems}, in
Mathematical {A}spects of the {F}inite {E}lement {M}ethod (I. Galligani, E. Magenes, eds), LNM vol. 606, Springer, New York, 1977, pp. 292--315. Zbl 0362.65089,
\bibitem{Sc} H. A. Schwarz, {Sur une d\'{e}finition errone\'{e} de l'aire d'une surface courbe}, in Gesammelte {M}athematische {A}bhandlungen, vol. 2, Springer, Berlin, 1890, pp.  309--311, 369--370.
\bibitem{Vo} M. Vohral{\'\i}k, {On the discrete Poincar\'{e}-Friedrichs inequalities for nonconforming approximations of the Sobolev space $H^1$}, Numer. Funct. Anal. Optim. \textbf{26:7-8} (2005), 925--952. Zbl 1089.65124,
\end{thebibliography}


\section{Appendix}\label{sec5}
\subsection{Proof of (\ref{Fried})}\label{sec51}
First of all, since nonconforming P1 element functions $w\in V_{\mathcal{T}}$ are piecewise linear, and can be parametrized by their edge midpoint values $w(M_e)$, $e\in \mathcal{E}$, we can explicitly estimate 
their discrete $H^1$ and $L_2$ norm:
\begin{equation}\label{VH1}
\sum_{\Delta\in \mathcal{T}} |\Delta| (\sum_{e \subset \Delta} |D_{e,\Delta}w|^2)\le 3\|w\|_{H^1,\mathcal{T}}^2,
\end{equation}
where the constant directional derivative $D_{e,\Delta}w$  of the linear function $w|_\Delta$ along the edge $e$ equals $2(w(M_{e'})-w(M_{e''}))/|e|$, where $e'$, $e''$ 
are the other two edges of $\Delta$. In the opposite direction, the inequality holds only with a constant depending on $\alpha_{\mathcal{T}}$.  Moreover,
\begin{equation}\label{VL2}
\|w\|_{L_2}^2 =\frac13 \sum_{\Delta\in \mathcal{T}} |\Delta| (\sum_{e \subset \Delta} |w(M_e)|^2).
\end{equation}

Consider all $2n+1$ triangles in the strip $\Omega_j=[0,1]\times 
[\frac{j-1}{2m},\frac{j}{2m}]$, and enumerate them consecutively starting from the left. Each $\Delta_i \in \Omega_j$, $i=0,\ldots,2n$, has exactly one edge (denoted $e_i$) parallel to the $x$-axis, and two edge midpoints (denoted by $M_i$ and $M_{i+1}$) on the line $y=\frac{2j-1}{4m}$. Obviously, for $i=1,\ldots,2n-1$,
we have 
$$
|w(M_{e_i})|\le \frac12 |w(M_{i+1})+w(M_i)| + |w(M_{e_i})-\frac12(w(M_{i+1})+w(M_i))|
$$
$$
\qquad \le \frac12( |w(M_{i+1})|+|w(M_i)|)+
\frac1{4m}|(w|_{\Delta_i})_y|,
$$
with the obvious modification 
$$
|w(M_{e_0})|\le |w(M_0)| + |w(M_{e_0})-w(M_0)| = |w(M_{0})| + \frac1{4m}|(w|_{\Delta_0})_y|,
$$
for $i=0$, and similarly for $i=2n$.
Thus, taking squares and using the inequality $(a+b)^2\le 2(a^2+b^2)$ the appropriate number of times, we get
$$
\sum_{i=0}^{2n} |\Delta_i| |w(M_{e_i})|^2 \le \frac{1}{2nm}\sum_{i=1}^{2n} |w(M_{i})|^2 + \frac1{8m^2}\sum_{i=0}^{2n} |\Delta_i| 
|(w|_{\Delta_i})_y|^2,
$$
and substitution gives
\begin{equation}\label{Fried1}
\frac13 \sum_{\Delta\subset \Omega_j} |\Delta| (\sum_{e \subset \Delta} |w(M_e)|^2)
\le \frac1{8m^2}\sum_{\Delta\subset \Omega_j} |w|^2_{1,\Delta} + \frac{1}{3nm} \sum_{i=1}^{2n} |w(M_{i})|^2.
\end{equation}

It remains to estimate the second term in (\ref{Fried1}). Since $w(M_0)=0$, we have
$$
\sum_{i=1}^n |w(M_{i})|^2=\sum_{i=1}^n \left|\sum_{l=1}^i (w(M_l)-w(M_{l-1})\right|^2 
\le \sum_{i=1}^n i\sum_{l=1}^i |w(M_l)-w(M_{l-1})|^2
$$
$$
\le \frac{n(n+1)}2\sum_{l=1}^n  |w(M_l)-w(M_{l-1})|^2.\qquad\qquad
$$
Now, take into account that
$$
|w(M_l)-w(M_{l-1})|^2=\frac1{4n^2}|(w|_{\Delta_{l-1}})_x|^2\le\frac{m}{n}|w|_{1,\Delta_{l-1}}^2,\qquad l=2,\ldots,n,
$$
and
$$
|w(M_1)-w(M_{0})|^2=\frac1{16n^2}|(w|_{\Delta_{l-1}})_x|^2\le\frac{m}{2n}|w|_{1,\Delta_{0}}^2,
$$
we see that
$$
\sum_{i=1}^{n} |w(M_{i})|^2 \le  \frac{m(n+1)}2\sum_{l=0}^{n-1} |w|_{1,\Delta_{l}}^2.
$$
In a similar fashion we also obtain
$$
\sum_{i=n+1}^{2n} |w(M_{i})|^2 \le  \frac{m(n+1)}2\sum_{l=n+1}^{2n} |w|_{1,\Delta_{l}}^2.
$$
Substitution into (\ref{Fried1}) gives
$$
\frac13 \sum_{\Delta\subset \Omega_j} |\Delta| (\sum_{e \subset \Delta} |w(M_e)|^2)
\le (\frac1{8m^2}+\frac{n+1}{6n})\sum_{\Delta\subset \Omega_j} |w|^2_{1,\Delta}<\frac12 \sum_{\Delta\subset \Omega_j} |w|^2_{1,\Delta}
$$
and, after summing up with respect to $\Omega_j$, $j=1,\ldots,2m$, according to (\ref{VH1}) and (\ref{VL2})  we arrive at (\ref{Fried}). 

\subsection{Proof of Lemma \ref{lem2}}\label{sec52}
We start with establishing (\ref{Wxy}) for all triangles interior to $\Omega'$ depicted in Figure \ref{fig2} on the left. 
By definition of the nodal values of $\tilde{w}$, and the fact that $u_{xy}$ is the product of two univariate linear polynomials, we compute
\begin{eqnarray*}
\tilde{w}_x^+&=&h^{-1}\left(\psi(x_0+\frac{h}{2},y_0+\frac{k}{2})-\psi(x_0-\frac{h}{2},y_0+\frac{k}{2})\right)\\
&=&-\left((2x_0+h)u_{xy}(x_0+\frac{h}{2},y_0+\frac{k}{2})+(2x_0-h)u_{xy}(x_0-\frac{h}{2},y_0+\frac{k}{2})\right)\\
&=& -4x u_{xy}(x_0,y_0+\frac{k}{2})-h\left(u_{xy}(x_0+\frac{h}{2},y_0+\frac{k}{2})-u_{xy}(x_0-\frac{h}{2},y_0+\frac{k}{2})\right)\\
&=& -4x_0u_{xy}(x_0,y_0+\frac{k}{2})-h^2u_{xxy}(x_0,y_0+\frac{k}{2})\\
&=& -(4x_0(1-2x_0)-2h^2)(1-2y_0-k),
\end{eqnarray*}
and, similarly,
$$
\tilde{w}_x^-=4x_0u_{xy}(x_0,y_0-\frac{k}{2})+h^2u_{xxy}(x_0,y_0-\frac{k}{2})=(4x_0(1-2x_0)-2h^2)(1-2y_0+k).
$$
Moreover,
\begin{eqnarray*}
\tilde{w}_y^+&=&k^{-1}\left(\psi(x_0+\frac{h}{2},y_0+\frac{k}{2})+\psi(x_0-\frac{h}{2},y_0+\frac{k}{2})\right)\\
&=&-hk^{-1} \left((2x_0+h)u_{xy}(x_0+\frac{h}{2},y_0+\frac{k}{2})-(2x_0-h)u_{xy}(x_0-\frac{h}{2},y_0+\frac{k}{2})\right)\\
&=& -hk^{-1}\left(2x_0(u_{xy}(x_0+\frac{h}{2},y_0+\frac{k}2)-u_{xy}(x_0-\frac{h}{2},y_0+\frac{k}2)) + 2h u_{xy}(x_0,y_0+\frac{k}2)\right) \\
&=& -2h^2k \left(u_{xy}(x_0,y_0+\frac{k}2)+x_0u_{xxy}(x_0,y_0+\frac{k}2)\right)=\kappa(4x_0-1)(1-2y_0-k),
\end{eqnarray*}
and
$$
\tilde{w}_y^- = -2h^2k \left(u_{xy}(x_0,y_0-\frac{k}2)+x_0u_{xxy}(x_0,y_0-\frac{k}2)\right) = \kappa(4x_0-1)(1-2y_0+k).
$$
This shows (\ref{Wxy}). The contribution of these triangles to the value of 
$|\tilde{w}|_{H^1,\mathcal{T}_{n,m}}^2$ (see (\ref{VH1}) for the formula) is of the order O$(1+\kappa^2)$.

For the triangles shown in Figure \ref{fig2} on the right, we have $\tilde{w}(M_e)=0$ for the horizontal and vertical edges,
which immediately leads to (\ref{WxyB}) if one substitutes the value for the remaining edge midpoint from (\ref{DefW}). 
This yields an  O$(h^2+\kappa^2)$ contribution to $|\tilde{w}|_{H^1,\mathcal{T}_{n,m}}^2$ from all triangles with sides
on the vertical boundaries of $\Omega$. 

It remains to check the triangles crossing the symmetry line $x=1/2$. Obviously, by the extension rule $\tilde{w}_x^\pm=0$ for all those triangles
while 
$$
\tilde{w}_y^{\pm} = \pm k^{-1}\psi(\frac12-\frac{h}2,y_0\pm \frac{k}2)= \pm 2\kappa(1-h)(1-2y_0\mp k).
$$
This gives (\ref{WxyM}). Consequently, we have to add another O$(\kappa^2)$ term to $|\tilde{w}|_{H^1,\mathcal{T}_{n,m}}^2$ which altogether
yields the desired estimate (\ref{WBound}) for the discrete $H^1$ norm of $\tilde{w}$.   Lemma \ref{lem2} is proved.

\subsection{Proof of Lemma \ref{lem3}}\label{sec53}
We give a bit more detail on the computations of the integrals involved than absolutely necessary.
For all triangles $\Delta^\pm$ but the ones depicted
in Figure \ref{fig2} on the right, in local coordinates, the linear function $\tilde{w}^{\pm}:=\tilde{w}|_{\Delta^\pm}$  equals 
$$
\tilde{w}^{\pm}(x_0+t,y_0+s)=\tilde{w}^\pm_x t + \tilde{w}^\pm_y s,\qquad  -h(1-k^{-1}|s|)\le t \le h(1-k^{-1}|s|),
$$
where $0\le s \le k$ for $\Delta^+$, and $-k\le s \le 0$ for $\Delta^-$. Therefore, we can use symmetries for triangle pairs $\Delta^\pm$ when evaluating
their contributions to $(f,\tilde{w})_{L_2}$. To do the calculations, we will use the following elementary formulas. For integers $\alpha,\beta\ge 0$ and the triangles $\Delta^\pm$ depicted in Figure \ref{fig2} on the left, we have
\begin{equation}\label{Beta}
I_{\alpha,\beta}^\pm:=\int_{\Delta^{\pm}_0} t^\alpha s^\beta \,dtds = \left\{\begin{array}{ll} 0,& \alpha \mbox{ odd},\\
&\\
(\pm 1)^\beta \frac{2 \alpha ! \beta !}{(\alpha+\beta+2)!}h^{\alpha+1}k^{\beta + 1},& \alpha \mbox{ even}.\end{array}\right.
\end{equation}
while for the triangles $\Delta^\pm$ depicted in Figure \ref{fig2} on the right it holds
\begin{equation}\label{BetaL}
\tilde{I}_{\alpha,\beta}^\pm:=\int_{\Delta^{\pm}_0} t^\alpha s^\beta \,dtds = 
(\pm 1)^\beta \frac{\alpha ! \beta !}{(\alpha+\beta+2)!}h^{\alpha+1}k^{\beta + 1}.
\end{equation}

Since, in local coordinates, 
$$
f(x_0+t,y_0+s)=2(x_0(1-x_0)+y_0(1-y_0) + (1-2x_0)t + (1-2y_0)s-t^2-s^2),
$$
using (\ref{Beta}) we compute
\begin{eqnarray*}
&&\int_{\Delta^{\pm}} f\tilde{w}\,dxdy = 
2\tilde{w}^\pm_x (1-2x_0) I_{2,0}^\pm\qquad\qquad\qquad\qquad\qquad\\
&&\qquad\quad  + 2\tilde{w}^\pm_y\left( ((x_0(1-x_0)+y_0(1-y_0))I_{0,1}+(1-2y_0)I_{0,2}^\pm-I_{0,3}^\pm-I_{2,1}^\pm\right)\\
&&\quad = \frac{h^3k}3(1-2x_0)\tilde{w}^\pm_x +\left(\pm \frac{2hk^2}3 (x_0(1-x_0)+y_0(1-y_0)) +\frac{hk^3}3 (1-2y_0) \mp(\frac{hk^4}{5}+\frac{h^3k^2}{15})\right)\tilde{w}^\pm_y.
\end{eqnarray*}
Thus,  the integral over $\Delta^{+}\cup \Delta^{-}$ equals
\begin{eqnarray*}
\int_{\Delta^{+}\cup \Delta^{-}} f\tilde{w}\,dxdy &=& hk\left(\frac{h^2}3(1-2x_0)(\tilde{w}^+_x+\tilde{w}^-_x) + \frac{k^2}3(1-2y_0)(\tilde{w}^+_y+\tilde{w}^-_y)\right.\\
&&\left. \quad  + \left(\frac{2k}3 (x_0(1-x_0)+y_0(1-y_0)) -\frac{1}{5}k^3-\frac{h^2k}{15}\right)(\tilde{w}^+_y-\tilde{w}^-_y)\right). \\
\end{eqnarray*}
Using (\ref{Wxy}) for $\tilde{w}^\pm_x$ and $\tilde{w}^\pm_y$, we have
$$
\tilde{w}^+_x+\tilde{w}^-_x=2(4x_0(1-2x_0)-2h^2)k=\mathrm{O}(k),
$$
and
$$
\tilde{w}^+_y+\tilde{w}^-_y=2\kappa(4x_0-1)(1-2y_0)=\mathrm{O}(h^2k^{-1}),\qquad
\tilde{w}^+_y-\tilde{w}^-_y=-2\kappa(4x_0-1)k=\mathrm{O}(h^2),
$$
and after substitution we see that each such integral is of order O$(h^3k^2)$. Consequently, the integral over the union of all such triangle pairs contained in $\Omega'$ is at most of order O$(h^2k)$. It is not hard to see that similar estimates hold 
for all triangles having an edge on one of the horizontal sides $y=0$, $y=1/2$ of $\Omega'$, or crossing the symmetry line $x=1/2$. 

For the triangles with $P=(0,y_0)$
depicted in Figure \ref{fig2} on the right, we have the following representations
in local coordinates: 
$$
f(t,y_0+s)=2(y_0(1-y_0) + t + (1-2y_0)s-t^2-s^2)
$$ 
and 
$$
\tilde{w}(t,y_0+s)=\tilde{w}^\pm + \tilde{w}^\pm_x t + \tilde{w}^\pm_y s,\qquad  0\le t \le h(1-k^{-1}|s|),
$$
where $0\le s \le k$ for $\Delta^+$, and $-k\le s \le 0$ for $\Delta^-$. Here, the absolute terms $\tilde{w}^\pm$ can be computed from the definition of $\tilde{w}$ as
$$
\tilde{w}^\pm=\pm h^2(1-h)(1-2y_0\mp  k)=\pm h^2(1-h)(1-2y_0) -h^2(1-h)k,
$$
while the derivatives $\psi^\pm_x$ and $\psi^\pm_y$ are given by (\ref{WxyB}).
Using (\ref{BetaL}) for the occuring integrals $\tilde{I}_{\alpha,\beta}^\pm$, we obtain
\begin{eqnarray*}
&&\int_{\Delta^{\pm}} f\tilde{w}\,dxdy = \tilde{w}_x^\pm \left(hky_0(1-y_0)+\frac{h^2k}3\pm
\frac{hk^2}3(1-2y_0)-\frac{h^3k+hk^3}6\right)\\
&&\qquad\qquad +\tilde{w}_x^\pm \left(\frac{h^2k}3y_0(1-y_0) + \frac{h^3k}6 \pm \frac{h^2k^2}{12}(1-2y_0)-\frac{3h^4k+h^2k^3}{30}\right)\\
&&\qquad\qquad + \tilde{w}_y^\pm \left(\pm\frac{hk^2}3 y_0(1-y_0) \pm \frac{h^2k^2}{12} + \frac{hk^3}{6}(1-2y_0)\mp\frac{h^3k^2+3hk^4}{30}\right).
\end{eqnarray*}
Due to (\ref{WxyB}) we have
$$
\tilde{w}^+_x+\tilde{w}^-_x=4hk(1-h)=\mathrm{O}(hk),\qquad 
\tilde{w}^+_x-\tilde{w}^-_x=4h(1-h)(1-2y_0)=\mathrm{O}(h),
$$
and
$$
\tilde{w}^+_y+\tilde{w}^-_y=-4\kappa(1-2y_0)(1-h)=\mathrm{O}(h^2k^{-1}),\qquad
\tilde{w}^+_y-\tilde{w}^-_y=4\kappa k(1-h)=\mathrm{O}(h^2).
$$
Together with the formula for $\tilde{w}^\pm$ this yields
\begin{eqnarray*}
\int_{\Delta^{+}\cup \Delta^{-}} f\tilde{w}\,dxdy 
&=& \frac{h^3k^2}3(1-2y_0)(1-h)(1-2y_0)\\
&& \quad  -\; h^3k^2(1-h)\left(y_0(1-y_0)+\frac{h}3-\frac{h^2+k^2}6\right)\\
&& \quad + \left(\frac{h^2k}3 y_0(1-y_0) + \frac{h^3k}6-\frac{3h^4k+hk^4}{30}\right)
(\tilde{w}^+_x+\tilde{w}^-_x) \\
&& \quad + \;\frac{h^2k^2}{12}(1-2y_0)(\tilde{w}^+_x-\tilde{w}^-_x) + \frac{hk^3}{6}(1-2y_0)(\tilde{w}^+_y+\tilde{w}^-_y)\\
&& \quad  + \left(\frac{hk^2}3(y_0(1-y_0) + \frac{h^2k^2}{12} -\frac1{30}(h^3k^2+3hk^4))(\tilde{w}^+_y-\tilde{w}^-_y)\right) \\
&& \quad \le C(h^3k^2 +h^3k^2 + h^2k^2 + h^3k^2+ h^3k^2)\le Ch^2k^2.
\end{eqnarray*} 
Summation with respect to all $\mathrm{O}(k^{-1})$ triangle pairs of this type gives another term of order O$(h^2k)$.
All in all we arrive at the statement of Lemma \ref{lem3}.

\end{document}